\documentclass[a4paper, 12pt]{article}

\usepackage{graphicx}
\usepackage{amssymb}

\def\blankbox{{ ~\hfill$\rlap{$\sqcap$}\sqcup$}}

\begin{document}

\bigskip \bigskip

\centerline {\bf Strongly Regular Graphs with No Triangles}

\bigskip \bigskip

\centerline{Norman Biggs}

\bigskip

\bigskip

\centerline{Department of Mathematics}

\centerline{London School of Economics}

\centerline{Houghton Street}

\centerline{London WC2A 2AE}

\centerline{U.K.}

\centerline{n.l.biggs@lse.ac.uk}

\bigskip

\centerline{Research Report  - September 2009}

\bigskip \bigskip

\centerline{\bf Abstract}
\bigskip

A  simplified version of the theory of strongly regular graphs is developed for the case
in which the graphs have no triangles.  This leads to (i) direct proofs of the Krein conditions,
and (ii) the characterization of strongly regular graphs with no triangles such that the second
subconstituent is also strongly regular. The method also provides an effective means of listing
feasible parameters for such graphs.
\vfill \eject

{\bf 1. Introduction}
\medskip

We consider graphs that are {\it strongly regular} and have {\it no triangles}, abbreviated to `SRNT graphs'.
Such a graph $X$ is characterized by two parameters $k$ and $c$, according to the rules
\smallskip

{\leftskip 30pt

$\bullet$  $X$ is regular with degree $k$;

$\bullet$  any two adjacent vertices have no common neighbours;

$\bullet$  any two non-adjacent vertices have $c$ common neighbours.

\par}

\smallskip

We shall discuss only cases where

$$k\ge 3, \qquad k > c \ge 1.$$

These conditions rule out the the pentagon ($k=2, c=1)$ and the
 complete bipartite graphs $K_{k,k}$, which have $c=k$.  Thus, by an SRNT graph we
mean a non-bipartite connected graph, with diameter $2$ and degree at least $3$.
Only six such graphs are known.
There are other pairs $(k,c)$ for which an SRNT graph may exist, but they are quite rare. Brouwer's list
[{\bf 2}] contains all possibilities with up to 280 vertices, and the Appendix to this paper
contains all possibilities with up to 6025 vertices.
\smallskip

If $v$ is a vertex of an SRNT graph $X$ we denote the subgraphs induced by the sets of vertices at distances
$1$ and $2$ from $v$ by $X_1(v)$ and $X_2(v)$ respectively.  Clearly, $X_1(v)$ has no edges, for any $v$.
It is worth remarking that the graphs $X_2(v)$ need not all be isomorphic, although we often use the
abbreviation $X_2$ to denote any one of them.
\smallskip

We shall need some standard notation and theory [{\bf 7}]. The sizes of $X_2$ and $X$ are given by
$$\ell =  \frac{k(k-1)}{c}, \qquad n = 1 + k + \ell = 1 + \frac{k}{c}(k-1+c).$$

The adjacency matrix $A$ of $X$ satisfies the equations $AJ = kJ$ and $A^2 + cA - (k-c)I = cJ$,
where $I$ is the identity matrix and
$J$ is the all-$1$ matrix. It follows that the eigenvalues of $A$ are
$k$ (with multiplicity 1) and the roots $\lambda_1, \lambda_2$
of the equation $\lambda^2 + c\lambda - (k-c) = 0$.  Furthermore, there is an integer $s >c$ such that
$c^2 + 4(k-c) =s^2$, where $s$ and $c$ have same parity, and the eigenvalues are the integers

$$ k = \frac {s^2 - c^2}{4} + c \qquad \lambda_1 = \frac{s-c}{2}, \qquad \lambda_2 = \frac{-s-c}{2}.$$

The multiplicities $m_1, m_2$ of $\lambda_1, \lambda_2$ are given by

$$m_1 = \frac{k}{2cs} \Big( (k-1+c)(s+c) - 2c\Big), \quad
  m_2 = \frac{k}{2cs} \Big( (k-1+c)(s-c) + 2c\Big).$$
Note that the algebra exhibits {\it s-symmetry}: replacing $s$ by $-s$ fixes $k$ but switches
$\lambda_1$ and $\lambda_2$, $m_1$ and $m_2$.
\smallskip

The conditions that $n$ and $m_1$ (and consequently $\ell$ and $m_2$) must be integers are known as
`feasibility conditions'.
They severely restrict the possible parameters $(k,c)$, for example, they imply that when $c \neq 2,4,6$
there are only finitely many possible pairs $(k,c)$.
\smallskip

In Section 3 we shall discuss another feasibility condition, which (among other things) greatly simplifies
the calculation of feasible parameters. If we take the basic parameters to be the positive integers $c$ and
$\lambda_1 =q$, then it follows that $s= c +2q$ and $k = (q+1)c + q^2$. The conditions that $n$ and $m_1$
are integers also take a fairly simple form in terms of $c$ and $q$ (details are given in the Appendix).
\smallskip

As well as the arithmetical conditions summarized above, it is possible to derive some `graph-theoretical'
conditions.  For example, a general result on distance-regular
graphs [{\bf 1}, {\bf 3}] implies that $k > 2c$, a result that can be proved in this case by simple arguments.
However, this result is not the best-possible. The feasibility condition given in Section 3
implies that $k \ge 3c-1$; indeed,
apart from the known examples, we must have $k \ge \frac{7}{2} c + \frac{25}{4}$.
\bigskip

{\bf 2. The second subconstituent  and its eigenspaces}
\medskip

It is clear from the definition of an SRNT graph $X$ that $X_2$ is a regular graph of degree $k-c$.
Simple arguments provide more specific information.
\medskip

{\bf Theorem 1} \quad $X_2$ is a connected graph with diameter  $2$ or $3$.
\smallskip

{\it Proof} \quad We show first that every path $uvw$ in $X$ is part of a $5$-cycle.
For any such path $v \in X_1(u)$ and $w \in X_2(u)$.  Let $x$ be any one of the $k-c$ vertices
in $X_2(u)$ that is adjacent to $w$, and let $y$ be any one of the $c$ vertices in $X_1(u)$  that is adjacent to
$x$.  Since $X$ has no triangles, $y \neq w$, and hence $uvwxy$ is a $5$-cycle.
\smallskip

Now let $d$ and $d_2$ denote the distance functions in $X$ and $X_2(u)$ respectively. If $p,q$
are vertices in $X_2(u)$ such that $d_2(p,q) >2$, then $d(p,q) = 2$. All $c$ vertices adjacent to
$p$ in $X_1(u)$ must also be adjacent to $q$, since $p$ and $q$ have $c$ common neighbours.
\smallskip

Any  path $qvp$ in $X$ with $v \in X_1(u)$ is part of a 5-cycle $qvpab$.
If $a$ were in $X_1(u)$ , then $a$ would also be adjacent to $q$ and we should have a triangle $abq$. Similarly
if $b$ were in $X_1(u)$ we should have a triangle $abp$.  Hence the path
$pabq$ must be in $X_2(u)$,  and $d_2(p,q) = 3$.
\blankbox
\medskip
\vfill \eject

In the case $c=1$ (known as the {\it Moore graph} case)
$X_2$ must have diameter $3$: in fact, $X_2$ is an antipodal $(k-1)$-fold covering of the complete graph
$K_k$ [{\bf 5}].  The known SRNT graphs with $k=3$ (Petersen) and $k=7$ (Hoffman-Singleton) exemplify
this situation.   For some other known SRNT graphs (Clebsch and Higman-Sims) $X_2$ has diameter $2$, and
indeed it is strongly regular and thus an SRNT graph.
\medskip

The algebraic theory of the subconstituents for any strongly regular graph is well-known [{\bf 7}, pp227-230].
In our case, when $X_1$ is trivial, it is possible to give a streamlined version.
\smallskip

Fix a vertex $v$ in $X$ and partition the adjacency matrix according to the vertex-partition
$\{v\} \cup X_1(v) \cup X_2(v)$:

$$A = \pmatrix{ 0 &J &O   \cr
                J &O   &B^T \cr
                O &B   &A_2 \cr }. $$
(Here, and in what follows, the $J$'s denote all-$1$ matrices of the appropriate sizes.)
\smallskip

For any $x \in {\mathbb R}^n$ let $[x_0 \; x_1 \; x_2]^T$ denote the corresponding column
vector, partitioned in the same way as $A$, so that $x_0 \in {\mathbb R}$, $x_1 \in {\mathbb R}^k$,
$x_2 \in {\mathbb R}^{\ell}$. By elementary matrix algebra it follows that
if $x$ is an eigenvector of $A$ with eigenvalue $\lambda \neq k$ then

$$Jx_1 = \lambda x_0, \quad x_0 J + B^T x_2 = \lambda x_1, \quad B x_1 + A_2 x_2 = \lambda x_2. $$

\medskip

{\bf Theorem 2} \quad  For any SRNT graph, $m_1 \ge k$ and $m_2 \ge k$.
\smallskip

{\it Proof} \quad Let $P$ be the eigenspace of $A$ for the eigenvector $\lambda_1$.
If $x \in P$ then $Jx = 0$, and since $Jx_1 = \lambda_1 x_0$, it follows that
$$0 = x_0 + Jx_1 + Jx_2 = (1 + \lambda_1)x_0 + J x_2.$$

Let $Q$ be the space comprising those $x$ for which
$x_0 = 0$ and $x_2 = 0$.  Since any $x \in Q$ also satisfies the equation $(1 + \lambda_1)x_0 + Jx_2 =0$,
it follows that the dimension of the space $P+Q$ is at most $n-1$.
\smallskip

On the other hand, if $x \in P \cap Q$, then $x_0 J + B^T x_2 = \lambda_1 x_1$, so $x=0$.
Thus, by a standard theorem, $\dim(P+ Q) = \dim P + \dim Q$, and

$$n-1 \ge \dim(P+Q) = \dim P + \dim Q = m_1 + k.$$

Since $n - 1 = m_1 + m_2$, it follows that $m_2 \ge k$.
\smallskip

Similarly, or by $s$-symmetry,  $m_1 \ge k$.
\blankbox
\medskip

Since $k-c$ is the degree of $X_2$, it is an eigenvalue of $A_2$, and since $X_2$ is connected, $k-c$
has multiplicity $1$.  We now determine the other possible eigenvalues of $X_2$.
Substituting the partitioned form of $A$ in the equation $A^2 + cA - (k-c)I = cJ$
we obtain three significant equations:

$$B^T B = (c-1)J + (k-1)I   \hskip 100pt (1)$$

$$ A_2^2 + cA_2 - (k-c)I + BB^T = cJ \hskip 70pt (2)$$

$$ A_2 B = -cB + cJ.  \hskip 148pt (3) $$

\medskip

{\bf Theorem 3} \quad Suppose $\mu \neq k-c$ is an eigenvalue of $X_2$.   Then either $\mu = \lambda_1$,
$\mu = \lambda_2$, or $\mu = -c$.
\smallskip

{\it Proof} \quad  If $x \neq 0$ is an eigenvector for $\mu \neq k-c$ then we have
$Jx = 0$.   By (2), $(A_2^2 + cA_2 - (k-c)I)x  = - BB^T x$. Since $A_2 x = \mu x$, we have
$$(\mu^2 +c \mu - (k-c))x = (\mu - \lambda_1)(\mu - \lambda_2)x = - BB^T x.$$
Thus if $B^T x = 0$, then $\mu = \lambda_1$ or $\mu = \lambda_2$.
\smallskip

Suppose $B^Tx  \neq 0$. Transposing (3) we have
$$B^T A_2 x =  -c B^T x, \quad {\rm that\; is} \quad \mu (B^T x) = -c (B^T x),$$
and so in this case $\mu = -c$.
\blankbox
\bigskip

{\bf 3. The Krein conditions and their consequences}
\medskip

Two more feasibility conditions involve parameters $K_1$ and $K_2$, known as the {\it Krein parameters}.
They arise in the general theory of distance-regular graphs, as described in [{\bf 3}]. After
some elementary algebra, in the strongly regular case they can be written in terms of
$k, \lambda_1, \lambda_2$ [{\bf 7}]:

$$K_1 = \lambda_1 \lambda_2^2 - 2\lambda_1^2 \lambda_2 - \lambda_1^2 - k\lambda_1 + k \lambda_2^2 + 2k\lambda_2
\hskip 50pt$$
$$      = (k + \lambda_1)(\lambda_2 + 1)^2 - (\lambda_1 + 1)(k + \lambda_1 + 2 \lambda_1 \lambda_2).$$

$$K_2 = \lambda_1^2 \lambda_2 - 2\lambda_1 \lambda_2^2 - \lambda_2^2 - k\lambda_2 + k \lambda_1^2 + 2k\lambda_1
\hskip 50pt$$
$$      = (k + \lambda_2)(\lambda_1 + 1)^2 - (\lambda_2 + 1)(k + \lambda_2 + 2 \lambda_1 \lambda_2).$$

When $X$ is a SRNT graph, we can  express $K_1$ and $K_2$ in terms of $s$ and $c$, as follows:

$$K_1 = \frac{1}{16}(s+c)(s-c+2)((s+c)^2 - 2(s+3c))$$
$$K_2 = \frac{1}{16}(s-c)(s+c-2)((s-c)^2 + 2(s-3c)).$$

Yet more elementary algebra leads to alternative, simpler, formulae.

$$   K_1 =  \frac{1}{4} (s+c)(s-c+2) (\lambda_2^2 + \lambda_2  -c) $$
$$       = \frac{1}{2} cs(s+c)\Big( \frac{m_1}{k} - 1 \Big).$$

$$   K_2 = \frac{1}{4} (s-c)(s+c-2)(\lambda_1^2 + \lambda_1  -c)$$.
$$    = \frac{1}{2} cs(s-c)\Big( \frac{m_2}{k} - 1 \Big).$$

These formulae provide a direct proof of the fundamental result on the Krein parameters, in the SRNT case.
\medskip

{\bf Theorem 4} \quad For any SRNT graph, $K_1 \ge 0$ with equality if and only if $m_1 = k$,
and $K_2 \ge 0$ with equality if and only if $m_2 =k$.
\smallskip

{\it Proof} \quad This  follows immediately from Theorem 2 and the formulae given above.
\blankbox

\bigskip

{\bf Corollary 1} \quad The matrix $A^2 + A - cI$ is positive semidefinite.
\smallskip

{\it Proof} \quad The formulae show that $\lambda_1$ and $\lambda_2$ both satisfy the condition
$\lambda^2 + \lambda - c \ge 0$, and the third eigenvalue $k$ also does so. (A direct proof of this
corollary may be possible.)
\blankbox

\medskip

{\bf Corollary 2} \quad  $k \ge 3c-1$.
\smallskip

{\it Proof} \quad We have shown that $\lambda_1^2 + \lambda_1  -c \ge 0$. Since
$\lambda_1$ satisfies the equation
$\lambda_1^2 + c\lambda_1 - (k-c) = 0$ it follows that
$$ -c \lambda_1 +(k-c) + \lambda_1 - c \ge 0, \quad {\rm that\;is,} \quad \lambda_1 \le \frac{k-2c}{c-1}.$$
Since $\lambda_1$ is a positive integer, $k-2c \ge c-1$.
\blankbox
\medskip

Parameters such as $(k,c) = (9, 4), (21, 10), \ldots\; $ are usually ruled out by calculating the
Krein parameters [{\bf 2}, {\bf 7}], but Corollary 2 achieves this result without any calculation.
Similar methods lead to the following general results.
\medskip

{\bf Corollary 3} \quad  The only SRNT graphs with $c+1 \le k \le 3c+4$ are the six currently-known ones.
\smallskip

{\it Proof} \quad We have to consider the cases $k = 3c + b$, $b = -1, 0, 1, 2,3,4$.
\smallskip

If $k = 3c -1$ we have $\ell = (3c-1)(3c-2)/c = 9c - 9 + 2/c$, so the only possiblities are $c=1$ and $c=2$.
These define the pentagon and the Clebsch graph.
\smallskip

Suppose $k= 3c + b$ with $b \ge 0$.  Then

$$s^2 = c^2 + 4(k-c) = c^2 + 8c + 4b = (c+4)^2 + 4(b-4).$$

When $b = 0,1,2,3$ this implies that $s^2 < (c+4)^2$, and since $s$ is a positive integer, $s \le c+3$.
Hence
$$(c+3)^2 \ge (c+4)^2 - 4(4-b), \quad {\rm that\;is} \quad c \le (9 - 4b)/2. $$

There are very few possibilities here, and the only one that gives a feasible set of parameters
is $b= 0, c=1, k=3$, which defines the Petersen graph.
\smallskip

If $k= 3c +4$ we have $\ell = (3c+4)(3c+3)/c = 9c - 21 + 12/c$, so $c$ is a divisor of $12$. The only feasible
solutions are $c=1,2, 4,6$, which define the graphs known by the names of Hoffman-Singleton, Gewirtz,
$M_{22}$, and Higman-Sims, respectively.
\blankbox

\medskip

{\bf Corollary 4} \quad An SRNT graph that is not currently-known must have $k \ge \frac{7}{2} c + \frac{25}{4}$.
\smallskip

{\it Proof} \quad In the light of the previous theorem, we can assume that $k = 3c+b$ with $b \ge 5$.
In this case $s^2 > (c+4)^2$ and hence $s \ge c+5$. Thus
$$(c+5)^2 \le c^2 +8c +4b, \quad {\rm that\;is} \quad c \le (4b - 25)/2. $$
In other words, $k = 3c + b$ with $b \ge \frac{1}{2}c + \frac{25}{4}$, as claimed.
\blankbox
\bigskip

{\bf 4. Linked pairs}
\medskip

We now consider {\it linked pairs} $(X,X')$ of SRNT graphs, that is, SRNT graphs
$X$ and $X'$ such that $X' = X_2(v)$ for every vertex $v$ of $X$. For comments on this problem,
see [{\bf 6}].
\medskip

{\bf Theorem 5} \quad The parameters of a linked pair $(X, X')$  are such that
\smallskip

$$k' = k-c, \qquad  c' = c-q, \qquad \lambda_1' = \lambda_1, $$
where
$$q = \frac{c^2(k-2)}{k^2 -(c+1)k + c(c-1)}.$$

{\it Proof} \quad  Clearly, the degree of $X'= X_2$ is $k' = k-c$. Since the number of
vertices of $X'$ is equal to $\ell$, we have
$$ 1 + k' + \frac{k'(k'-1)}{c'} = \frac{k(k-1)}{c}.$$
Substituting $k' = k-c$ and solving for $c'$ gives $c' = c-q$, where $q$ is as stated.
\smallskip

Since $X'$ is strongly regular, it has two eigenvalues other than $k-c$, and just one of
them  ($\lambda_1'$) is positive.  According to Theorem 3, the only possible eigenvalues are $\lambda_1$,
$\lambda_2$ and $-c$, of which only $\lambda_1$ is positive. Hence $\lambda_1' = \lambda_1$.
\blankbox
\medskip

Combining these equations leads to our main result.  Since $\lambda_1' = \lambda_1$ it follows that
$s'-c' = s - c$, and hence $s' = s-q$. From the equations
$s'^2 = c'^2 + 4(k'- c')$ and $s^2 = c^2 + 4(k-c)$  we obtain
$$ (s-q)^2 = (c-q)^2 + 4(k-c -c +q) = s^2 - 2qc + +q^2 -4c +4q,$$
$$ {\rm that\; is} \quad s = c-2 + \frac{2c}{q}.$$

Here both $s$ and $q$ can be written as functions of
$k$ and $c$. This yields the equation
$$c^2 (k-2)^2 (c^2 + 4k - 4c) = \Big( c(c-2)(k-2) + 2k^2 - 2(c+1)k +2c(c-1) \Big) ^2,$$

which is a quartic in $k$, and factors conveniently:

$$4(k-1)(k-c) (k^2 - (3c+1)k  - c(c^2 -4c -1)) = 0.$$

Thus, if there is a linked pair $(X,X')$ and $c$ is given, $k$ must be
a positive integer root of the quadratic factor. The discriminant of this factor is

$$\Delta = (3c+1)^2 + 4c(c^2 - 4c -1) = (c-1)^2 (4c + 1).$$

So $\Delta$ is a perfect square if and only if $4c+1$ is the square of an integer, which must be an odd
number $2r+1$. That is, $c= r(r+1)$.
The corresponding value of $k$ is

$$ \frac{1}{2} ( 3c +1 + \sqrt{\Delta}) \; = \;  r(r^2 + 3r + 1),$$
and these are the only values for which a linked pair can exist. Furthermore, for these values

$$q \; = \; \frac{c^2(k-2)}{k^2 -(c+1)k + c(c-1)} \; = \; r.$$

It is easy to check that $k = q(q^2 + 3q +1)$ and $c = q(q+1)$ satisfy all the feasibility
conditions for an SRNT graph $X$, as do the corresponding values for $X'$, $k' = q^2(q+2)$ and
$c' = q^2$.  Precisely, we have

$$\ell = (q^2 + 2q -1)(q^2 + 3q +1), \quad  n = q^2(q+3)^2, \quad  s = q(q+3), $$
$$\lambda_1 = q, \quad \lambda_2 = -q(q+2) , \quad m_1 = (q^2 +2q -1)(q^2 +3q +1),$$
$$ m_2 = q(q^2 + 3q + 1), \quad K_1 = q^2(q+1)(q+2)(q+3)(q^2 +q -1),\qquad  K_2 = 0.$$
$$\ell' = (q+1)(q+2)(q^2 + q - 1), \quad  n' = (q^2 + 2q -1)(q^2 + 3q +1) , \quad  s' = q(q+2), $$
$$\lambda_1' = q, \quad \lambda_2' = -q(q+1) , \quad m_1' = (q^2 + 3q + 1)(q^2 + q -1),$$
$$ m_2' =   (q+1)(q^2 +2q -1), \quad K_1' = q^2(q+1)^2(q^3 + 2q^2 -q-1), \quad K_2' = q^2(q^2 + q -1).$$

Thus we have the main result.
\medskip

{\bf Theorem 6} \quad The parameters of a linked pair $(X,X')$ of SRNT graphs must be of the form
$$ k = q(q^2 + 3q + 1), \quad c = q(q+1), \quad k' = q^2(q+2), \quad c' = q^2,$$
where $q$ is a positive integer.  Both sets of parameters are feasible for all $q \ge 1$.
(When $q=1$ we obtain the Clebsch/Petersen pair, and when $q=2$ we obtain the Higman-Sims/$M_{22}$ pair.
These graphs are known to be the unique ones with the relevant parameters.)
\blankbox
\medskip

Similar results have been obtained by Cameron [{\bf 4}, Theorem 5] and Smith [{\bf 11}, Theorem E].
Cameron used a result on partial quadrangles, and Smith considered the
case when $X$ admits a group of automorphisms that acts transitively
on the vertices, and the
stabilizer of a vertex $v$ acts transitively as a group of automorphisms of $X' = X_2(v)$. Her proof involves
calculations with the constituents of the
permutation characters, which appear similar to the calculations given above.
\smallskip

 The values for $k$ and $c$ are the SRNT case of the family
known as {\it negative latin square} parameters, first obtained by Mesner [{\bf 8}]. Graphs of this type
were also studied by M. Shrikhande [{\bf 9}] and S. Shrikhande [{\bf 10}].

\vfill \eject

{\bf References}
\medskip

{\bf 1.} N.L. Biggs.  Automorphic graphs and the Krein condition. {\it Geom. Dedicata} (5) 1976 117-127.

{\bf 2.} A.E. Brouwer.  Strongly Regular Graphs. In: {\it Handbook of Combinatorial Designs}, ed. C. Colbourn,
J. Dinitz,  CRC Press, 1996.

{\bf 3.} A.E. Brouwer, A.M. Cohen, A.Neumaier.  {\it Distance-Regular Graphs}, Springer, Berlin 1989.

{\bf 4.} P.J. Cameron. Partial quadrangles. {\it Quart. J. Math. Oxford (2)} 26 (1975) 61-73.

{\bf 5} A.D. Gardiner. Antipodal covering graphs {\it J. Combinatorial Theory (Series B)} 16 (1974) 255-273.

{\bf 6.} C.D. Godsil. Problems in algebraic combinatorics. {\it Elect. J. Combinatorics} 2 (1995) F1.

{\bf 7.} C.D. Godsil, G. Royle. {\it Algebraic Graph Theory}, Springer, New York 2001.

{\bf 8.} D.M. Mesner. A new family of partially balanced incomplete block designs with some latin square
design properties. {\it Ann. Math. Statist.} 38 (1967) 571-581.

{\bf 9.} M.S. Shrikhande. Strongly regular graphs and quasi-symmetric designs. {\it Utilitas Mathematica}
3 (1973) 297-309.

{\bf 10.} S.S. Shrikhande. Strongly regular graphs containing strongly regular subgraphs.
{\it Proc. Indian Natl. Sci. Acad. Part A} 41 (1975) 195-203.

{\bf 11.} M.S.Smith. On rank 3 permutation groups.  {\it J. Algebra} 33 (1975) 22-42.

\bigskip

{\bf Appendix}
\medskip

This Appendix contains calculations from Nimashi Thilakaratne's  dissertation (2009) for the
MSc in Applicable Mathematics at the LSE.  The calculations are based on the following result (see also
Cameron [{\bf 4}]).
\medskip

{\bf Theorem} \quad The number $n$ of vertices of an SRNT graph with $\lambda_1 = q$ is in the range

$$ \Bigl\lceil   2q^3 + 3q^2 - q +2q(q+1) \sqrt{q^2 +q -2} \Bigr\rceil \;  \le  \; n \;  \le \;  q^2 (q+3)^2.$$

The parameter $c$ is in the range $1 \le c \le q(q+1)$, and must be such that
\smallskip

{\leftskip 30pt
$c$ is a divisor of $q^4 - q^2$, and
\smallskip

$c + 2q$ is a divisor of  $q^4 +3q^3 + 5q^2 +3q  +  q(q^4 - q^2)/c$.
\par}

\medskip

{\it Proof} \quad  Given the values of $\lambda_1 = q$ and $c$, the parameters $s$, $k$ and $n$ are

$$s = c+2q, \qquad k = (q+1)c + q^2, \qquad n = Ac + B + D/c,$$

where $A = q^2 + 3q +2$, $B = 2q^3 + 3q^2 - q$, $D = q^4 - q^2$.
\smallskip

According to the formulae given in Section 3, the condition $K_2 \ge 0$ implies that
$\lambda_1^2 + \lambda_1 - c \ge 0$.  Hence
$c$ lies in the range $1 \le c \le  q(q+1)$. As a function of $c$,
$n$ has only one extreme point, a minimum, at the point where
$$c^2 = D/A, \quad  {\rm that\; is} \quad  c = q \left(\frac{q-1}{q+2} \right) ^{\frac{1}{2}}.$$
Substituting this value of $c$ gives the minimum value $n_{min}$, and since $n$ must be an integer,
we get the result as stated above.
\smallskip

The maximum value $n_{max}$ must therefore occur at one of the ends of the range, and
a simple calculation shows that the values at $c=1$ and $c=q(q+1)$ respectively are
$$ q^4 +2q^3 + 3q^2 + 2q +2 \qquad {\rm and} \qquad q^2(q+3)^2.$$
So the maximum occurs when $c=q(q+1)$.
\smallskip

If such a graph exists, $n$ and $m_1$ must be integral.  Another calculation gives
$$ m_1 = Ac + E + \frac{Fc + qD}{c(c + 2q)} ,$$
where $E = q^3 - 4q -2$, $F= q(q+1)(q^2 + 2q + 3)$, and $A, D$ are as above.
\smallskip

If $n$ is integral $c$ must divide $D$.  In that case $c$ must also divide $Fc + qD$ and hence
the condition that $m_1$ is an integer reduces to the fact that $c+2q$ must divide $F + q(D/c)$.
\blankbox
\bigskip

This theorem enables feasible parameters to be calculated systematically. The method is to
fix $q$ and find those $c$ in the range $1 \le c \le q(q+1)$ such that $c$ and $c+2q$ satisfy the divisiblity
conditions. For example, when $q=4$ we require the integers $c$ such that $1 \le c \le 20$,
 $c$ divides $240$, and $c+8$ divides $540 + 960/c$.  It is easy to check that the
 only possibilities are $c= 2,4,6,12,16,20$.
\smallskip

The theorem also gives bounds $n_{min}$ and $n_{max}$, and these provide an effective method of
tabulating the results.  For $1 \le q \le 11$ the bounds are as follows:

$$ \matrix{    q         &1 &2  &3   &4   &5   &6   &7    &8    &9    &10   &11    \cr
            n_{min}   &4 &50 &154 &342 &638 &1066&1650 &2413 &3381 &4577 &6025  \cr
            n_{max}   &16&100&324 &784 &1600&2916&4900 &7744 &11664&16900&23716 \cr}. $$

Suppose we wish to list all the feasible parameters for SRNT graphs with at most 1000 vertices.
According to the table, we need only carry out the calculation for $1 \le q \le 5$, since
$n_{min}(6)$ is greater than 1000. Similarly, if we list the feasible parameters for $1 \le q \le 10$,
the list will contain all possibilities with fewer than 6025 vertices. The results of these calculations
are tabulated below.

$$\matrix{ n &k &c &s &\ell &\lambda_1 & \lambda_2 &m_1 &m_2 &K_1 &K_2 \cr
             &&&&&&&&&& \cr
           10  &3   &1  &3 &6    &1         &-2         &5   &4   &4   &1   \cr
           16  &5   &2  &4 &10   &1         &-3         &10  &5   &24  &0   \cr
           50  &7   &1  &5 &42   &2         &-3         &28  &21  &45  &20   \cr
           56  &10  &2  &6 &45   &2         &-4         &35  &20  &120 &24   \cr
           77  &16  &4  &8 &60   &2         &-6         &55  &21  &468 &20   \cr
           100 &22  &6  &10&77   &2         &-8         &77  &22  &1200 &0   \cr
           162 &21  &3  &9  &140     &3          &-6           &105    &56    &648    &135   \cr
           176 &25  &4  &10 &150     &3          &-7           &120    &55    &1064   &144   \cr
           210 &33  &6  &12 &176     &3          &-9           &154    &55    &276    &144   \cr
           266 &45  &9  &15 &220     &3          &-12          &209    &56    &5904    &99   \cr
           324 &57  &12 &18 &266     &3          &-15          &266    &57    &11880    &0   \cr
           352 &26  &2  &10 &325     &4          &-6           &208    &143   &840    &360   \cr
           352 &36  &4  &12 &315     &4          &-8           &231    &120   &2080   &448   \cr
           392 &46  &6  &14 &345     &4          &-10          &276    &115   &4200   &504   \cr
           552 &76  &12 &20 &475     &4          &-16          &437    &114   &18240  &480   \cr
           638 &49  &4  &14 &588     &5          &-9           &406    &231   &3672  &1040   \cr
           650 &55  &5  &15 &594     &5          &-10          &429    &220   &5100  &1125   \cr
           667 &96  &16 &24 &570     &4          &-20          &551    &115   &36400  &304   \cr
           704 &37  &2  &12 &666     &5          &-7           &407    &296   &1680   &840   \cr
           784 &116 &20 &28 &667     &4          &-24          &667    &116   &63840    &0   \cr
           800 &85  &10 &20 &714     &5          &-15          &595    &204   &18000 &1400   \cr}$$

\bigskip
\centerline{\it Table 1:  Feasible parameters for SRNT graphs with at most 1000 vertices}

\vfill \eject

\centerline{$q=5$ (including those listed in Table 1)}

$$\matrix{ c &2   &4   &5   &10   &20   &25  &30   \cr
           k &37  &49  &55  &85   &145  &175 &205  \cr
           n &704 &638 &650 &800  &1190 &1394&1600 \cr } $$
\medskip

\centerline{$q=6$}

$$\matrix{ c &2    &4    &6    &9    &15   &30   &36  &42   \cr
           k &50   &64   &78   &99   &141  &246  &288    &330    \cr
           n &1276 &1073 &1080 &1178 &1458 &2256 &2585&2916     \cr } $$

\medskip

\centerline{$q=7$}

$$\matrix{ c &1    &4   &6     &7    &14   &21    &28   &42   &49  &56   \cr
           k &57  &81   &97    &105  &161  &217   &273  &385  &441 &497  \cr
           n &3250&1702 &1650  &1666 &2002 &2450  &2926 &3906 &4402&4900     \cr } $$

\medskip

\centerline{$q=8$}

$$\matrix{ c &2    &4    &6     &8     &14   &24  &28  &56  &64  &72   \cr
           k &82   &100  &118   &136   &190  &280 &316 &568 &640 &712    \cr
           n &3404 &2576 &2420  &2432  &2756 &3536&3872&6320&7031&7744    \cr } $$
\medskip

\centerline{$q=9$}

$$\matrix{ c &2    &4    &9    &12   &15   &18   &27   &36   &72   &81     &90  \cr
           k &101  &121  &171  &201  &231  &261  &351  &441  &801  &891    &981   \cr
           n &5152 &3752 &3402 &3552 &3774 &4032 &4902 &5832 &9702 &10682  &11664  \cr } $$

\medskip

\centerline{$q=10$}

$$\matrix{ c &2    &4    &6    &10    &20   &45  &90    &100   &110    \cr
           k &122  &144  &166  &210   &320  &595 &1090  &1200  &1310   \cr
           n &7504 &5293 &4732 &4600 &5425 &8450 &14280 &15589 &16900    \cr } $$

\bigskip

\centerline{\it Table 2:  Feasible parameters for SRNT graphs with $\lambda_1 = q = 5,6,7,8,9,10$}

\end{document}